\documentclass[11pt]{article}\include{ddefs}
\usepackage{graphics,pifont,array}
\usepackage{graphicx}
\usepackage{endnotes}
\let\footnote=\endnote

\begin{document}

\begin{titlepage}

\begin{center}

{\Large{        \mbox{   }                          \\
                \mbox{   }                          \\
                \mbox{   }                          \\
                \mbox{   }                          \\
                \mbox{   }                          \\
                \mbox{   }                          \\
                \mbox{   }                          \\
                \mbox{   }                          \\
                \mbox{   }                          \\
                \mbox{   }                          \\
                \mbox{   }                          \\
  {\textbf{RADIATIVE ALBEDO FROM A                  \\
          LINEARLY FIBERED HALF SPACE}}             \\
               \mbox{    }                          \\
               \mbox{    }                          \\
              J. A. Grzesik                         \\
           Allwave Corporation                      \\
        3860 Del Amo Boulevard                      \\
                Suite 404                           \\
           Torrance, CA 90503                       \\
                \mbox{    }                         \\
           (818) 749-3602                           \\ 
            jan.grzesik@hotmail.com                 \\
              \mbox{     }                          \\
              \mbox{     }                          \\  
              \today                                      }  }

\end{center}

\end{titlepage}

\setcounter{page}{2}

\pagenumbering{roman}
\setcounter{page}{2}
\vspace*{+2.725in}

\begin{abstract}

\parindent=0.245in

    A growing acceptance of fiber reinforced composite materials imparts some relevance
    to exploring the effects which a predominantly linear scattering lattice may have
    upon interior radiative transport.  Indeed, a central feature of electromagnetic
    wave propagation within such a lattice, if sufficiently dilute, is ray confinement
    to cones whose half-angles are set by that between lattice and the incident ray.
    When such propagation is subordinated to a viewpoint of photon transport, one arrives
    at a somewhat simplified variant of the Boltzmann equation with spherical scattering
    demoted to its cylindrical counterpart.  With a view to initiating a hopefully wider
    discussion of such phenomena, we follow through in detail the half-space albedo
    problem.  This is done first along canonical lines that harness the Wiener-Hopf
    technique, and then once more in a discrete ordinates setting via flux decomposition
    along the eigenbasis of the underlying attenuation/scattering matrix.  Good agreement
    is seen to prevail.  We further suggest that the Case singular eigenfunction apparatus
    could likewise be evolved here in close analogy to its original, spherical scattering
    model.

\end{abstract}

\pagestyle{plain}

\parindent=0.5in

\newpage

\pagenumbering{arabic}

\pagestyle{myheadings}

\setlength{\parindent}{0pt}

\pagestyle{plain}

\parindent=0.5in

\newpage
\mbox{   }

\pagestyle{myheadings}

\markright{J. A. Grzesik \\ radiative albedo from a linearly fibered half space}

\section{Introduction}

     Fiber reinforced composites\footnote{Abbreviated henceforth as FRC, singular or
plural as the local context may dictate.} have, in recent decades, enjoyed a widespread
penetration into the manufacture of durable structures, most noticeable among them being
perhaps aircraft fuselage segments on a large scale, an arena dominated heretofore by
aluminum as the material of choice.  And, while the primary objectives of reduced weight
and cost have indeed been achieved, subsidiary issues of somewhat lesser importance have
nevertheless arisen.  One among the latter is FRC behavior {\textit{vis-{\`{a}}-vis}
radiative transport at $\mu$m wavelengths in the infrared.

      Reinforcing fibers within some uniform host matrix are most frequently encountered
at two opposite extremes of dispersal regularity/irregularity, to wit, short tendrils
randomly oriented, or else extended threads aligned along a single direction.  Evidently
it is only the FRC samples which belong to the second, aligned category which can be
expected to exert any significant orientational influence upon radiative transport.

      With this background in mind, the present note seeks to offer a very modest first
insight into FRC radiative transfer by solving the albedo problem for an otherwise
uniform half-space matrix\footnote{Admittedly, in what follows, we allow this matrix
to default to empty space.  This is done so as to avoid having interface refractive complications
obscure the details of our evolving radiative transport machinery.  Needless to say,
this physical/mathematical lacuna pleads to be filled by more realistic work in the future.}
randomly seeded with long (idealized, to be sure, as infinitely long)
fibers in sufficiently dilute\footnote{By dilute we mean here an average fiber
separation of at least several wavelengths.} distribution parallel to the material
interface.

     It is characteristic of an electromagnetic field impinging upon a uniaxial scatterer
as now described to propagate along directions confined to cones having the fiber direction
as axis and opened to the angle between that axis and the direction of incidence.\footnote{This
geometrically engaging attribute is, regrettably, passed over in silence by the standard
electromagnetics texts (such as those authored by Stratton, Jackson, Panofsky and Phillips, and
so on).  Sheltered though it may thus be, its persistent presence in the literature is securely
anchored around [{\textbf{1}}] and in research papers of the type [{\textbf{2-4}}] to which it
lends support.  One finds there a derivation which relies upon the peculiarities of an
essentially two-dimensional field propagation, complete with reference to the asymptotic
behavior of the Hankel functions which figure therein.  By contrast, an {\textit{ab initio}}
derivation from first principles, altogether liberated from any debt to special function features,
can be found in [{\textbf{5}}].}  Since conical confinement of this sort is indifferent to the
incident field polarization, a welcome invitation presents itself to sidestep the Maxwell
electromagnetic apparatus in favor of the much simpler Boltzmann equation governing radiative
transfer.  Such transfer need clearly be tracked only in its projection upon the plane perpendicular
to fiber direction, the full conical flux being gotten therefrom under a simple multiplication
by $\sec\phi$ as indicated in Figure 1.\footnote{Figure 1 is freely available from the internet
under a Google search ``cylinder cone scattering."  Its origin can be traced to [{\textbf{6}}], and it
reappears on p. 264 of [{\textbf{2}}].}  Plainly put, the only positional variables
with which the Boltzmann transport equation need concern itself are the angle $\theta$ about fiber
direction (Figure 1) and the depth of penetration, here taken as co\"{o}rdinate $y$, into the
fiber laden half space.

      Our objectives here are very limited indeed, barely embryonic if one be permitted to say so.
On the one hand we shall not venture into the prevailing technology of FRC manufacture, while, on the other,
our analysis will lay claim merely to a methodological preamble to more refined, more sophisticated \mbox{compu-}
tations which may yet appear in the future from
the pens of other hands.  Indeed, on the methodological side we shall restrict ourselves, as
already stated, to a uniform, uniaxially fibered half space, and, moreover, to a scattering
function within Boltzmann's equation which is invariant against azimuthal angle $\theta$ around fiber
direction (co\"{o}rdinate axis $z$ as indicated in Figure 1).  Our primitive goal will
be the interface albedo under plane wave irradiation at azimuthal angle $\theta_{0}$
(and its cosine $\mu_{0}=\cos\theta_{0}$), reckoned from the positive $y$ axis and
arbitrarily chosen in the impact range $-\pi/2<\theta_{0}<\pi/2.$\footnote{In keeping with the
discussion above, the value of oblique incidence angle $-\pi/2<\phi_{0}<\pi/2$ from Figure 1 should
of course be retained in the recesses of one's mind, but need not otherwise be allowed to enter the
details of Boltzmann equation calculations.}

       Within these strict confines we shall pursue our albedo target along two routes, one
quasi-numerical, based upon the Boltzmann equation being discretized along angular nodes
conforming to a Gaussian quad-
\newpage
\mbox{   }
\mbox{   }
\newline
\newline
\newline
rature net, the other more traditionally analytic and relying
upon a Wiener-Hopf complex plane machinery.  Essentially perfect agreement will be demonstrated
between the albedo outcomes at the end of these distinct routes, and brought to the fore will
also be a standard incoming/outgoing symmetry that harkens back to the Ambarzumian/Chandrasekhar/Case
machinery of invariant embedding [{\textbf{7-10}}].  All of this, it is hoped, may yet stimulate
more refined, more realistic analyses by other investigators engaged in radiant transport across fibered media.
\vspace{-0.12in}
\begin{center}
\includegraphics[width=0.45\linewidth]{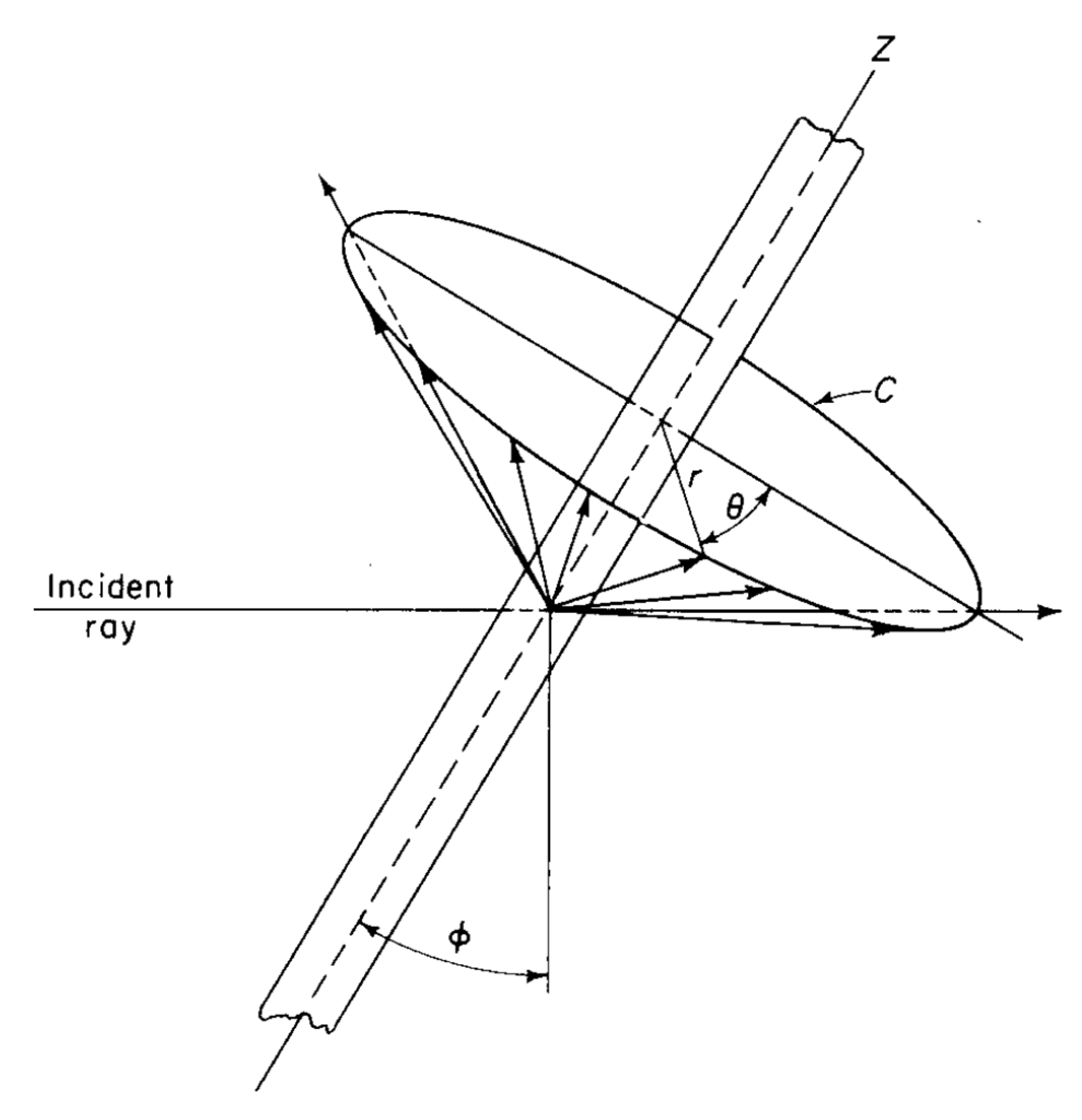}
\end{center}
\vspace{-15mm}
\begin{center}
{\large{
Figure 1.  Cylinder scattering cone C under oblique ray impact } }
\end{center}
\parindent=0.5in

      The material thus described is organized as follows.  After an obligatory setup of albedo geometry and the
      underlying Boltzmann equation, pride of place is given to the Wiener-Hopf machinery.  Sketched next are the
      very simple details of transport equation discretization along the nodes of a concatenated Gaussian quadrature
      lattice extended across the full, $-\pi<\theta<\pi$ range of scattering azimuth $\theta.$  \mbox{There fol-}      
      low polar albedo plots across the outgoing range $\pi/2<|\theta|<\pi,$ with special emphasis
      upon a virtually perfect agreement between calculations having an entirely dissimilar provenance.  Our note
      then concludes with one or two suggestions as to possible extensions of this work.

\section{Albedo Setup}

 Right-handed Cartesian axes, their origin in the half-space boundary, are placed so that axis $z$ is aligned
 with fibers as in Figure 1, axis $y$ points into the scattering medium, while axis $x$ is regarded as the
 vertical.  As already indicated, azimuthal angle $\theta,$ $-\pi<\theta<\pi,$ is reckoned from axis $y,$
 positive or negative when conveying a right-hand rotation around axis $\mp z.$
 
  \newpage
 \mbox{   }
 \newline

       Powerless to change, we betray an atavistic fidelity to the nomenclature of neutron transport by
       writing now $\psi(\theta,y)$ for the photon flux per unit increment of scattering angle $\theta.$  The total
       macroscopic cross section (absorption plus scattering, $\sigma_{a}+\sigma_{s}$;  dimension of an inverse length)
       is written as $\sigma,$ with both of its constituents assumed to be spatially constant as to both position $y$ and
       angle $\theta.$  Normalization by $\sigma$ yields the dimensionless single-scattering albedo $\omega=\sigma_{s}/
       \sigma,$ $0\leq \omega <  1,$\footnote{$\rule{0mm}{0.75in}$We shall expressly assume
       here that the fibers absorb to a greater or lesser
       extent, so that indeed $\omega<1.$  Calculation of the macroscopic parameters $\sigma_{a,s}$ is a considerable
       electromagnetic/statistical task in its own right, one that we do not propose to undertake here.
       A partial hint, but only a hint, as to how such onerous computations must proceed can be found in [{\bf{11}}]
       and the references there cited.}  and an optical path length $\tau=\sigma y$ gotten when $y$ is measured in units of one mean free path.
       Playing a key r\^{o}le is a probability $P(\theta|\theta')/(2\pi)d\theta$ that a photon arriving along
       direction $\theta'$ is scattered into angular increment $d\theta$ around $\theta.$  One evidently must insist
       upon the normalization $\int_{-\pi}^{\,\pi}P(\theta|\theta')d\theta=2\pi,$ regardless of incidence angle
       $\theta'.$  In what follows we simply default to the primitive, isotropic case $P(\theta|\theta')=1.$  In
       standard notation, $\mu=\cos\theta.$
       
       
       
       With these obligatory definitions duly disposed of, we may simply set down
       \begin{equation}
       \mu\frac{\partial \psi(\theta,\tau)}{\partial \tau}+\psi(\theta,\tau)=\frac{\omega}{2\pi}
       \int_{-\pi}^{\,\pi}\psi(\theta',\tau)d\theta'
       \end{equation}
       as our basic photon bookkeeping (Boltzmann) equation in the $x-y$ plane transverse to the mandated fiber direction.
       Flux $\psi(\theta,\tau)$ is further decomposed as the sum of a Dirac delta function sheet
       $e^{-\tau/\mu_{0}}\delta(\theta-\theta_{0})$ incident at angle $\theta_{0}$ and properly attenuated thereafter,
       and a diffuse, multiply scattered component $\psi_{d}(\theta,\tau).$  Thus
       \begin{equation}
       \psi(\theta,\tau)=e^{-\tau/\mu_{0}}\delta(\theta-\theta_{0})+\psi_{d}(\theta,\tau)\,.
       \end{equation}
       Decomposition (2) converts the homogeneous equation (1) into an inhomogeneous counterpart
       \begin{equation}
       \mu\frac{\partial \psi_{d}(\theta,\tau)}{\partial \tau}+\psi_{d}(\theta,\tau)=\frac{\omega}{2\pi}
       \int_{-\pi}^{\,\pi}\psi_{d}(\theta',\tau)d\theta'+\frac{\omega}{2\pi}e^{-\tau/\mu_{0}}
       \end{equation}
       whose source $(\omega/2\pi)e^{-\tau/\mu_{0}}$ is now entirely overt.  And then, as a boundary condition on (3)
       we further require that there be no re\"{e}trant diffuse flux at half-space boundary $\tau=0,$ {\em{viz.,}}
       $\psi_{d}(\theta,0)=0$ when $|\theta|\leq \pi/2.$  By contrast, it is the structure of the escaping, diffuse flux
       $\psi_{d}(\theta,0)$ for $\pi/2<|\theta|<\pi$ that we identify with the half-space albedo.

       \section{Wiener-Hopf solution for \mbox{\boldmath{$\psi_{d}(\theta,\tau)$}}}
       
             We displace attention from $\psi_{d}(\theta,\tau)$ {\em{per se}} to its angular aggregate
             \[\rho_{d}(\tau)=\int_{-\pi}^{\,\pi}\psi_{d}(\theta,\tau)\,d\theta\]
             which, following division by the speed of light $c,$ simply measures the diffuse photon
             density.  It turns out that nothing is lost thereby since, on the basis of (3), we can
             further write
             \begin{equation}          
             \psi_{d}(\theta,\tau)=-\frac{\omega}{2\pi\mu}\int_{\tau}^{\,\infty}\rho_{d}(\theta,\tau')
             e^{(\tau'-\tau)/\mu}d\tau'+\frac{\omega}{2\pi}\frac{\mu_{0}}{\mu_{0}-\mu}e^{-\tau/\mu_{0}}
             \end{equation}
             throughout the entire retrograde flight range $\pi/2<|\theta|<\pi$ wherein $-1<\mu<0.$  On setting
             $\tau=0$ we see that knowledge about $\rho_{d}(\tau)$ provides a direct stepping stone to the desired
             albedo.
             \newpage
             \mbox{    }
             \newline 
             
             A path to solution for both $\psi_{d}(\theta,\tau)$ and $\rho_{d}(\tau)$ is cleared by first subjecting
             these quantities to Laplace transformation, denoted by a circumflex, {\em{viz.,}}
             \begin{equation}
             \hat{\psi}_{d}(\theta,s) = \int_{0}^{\,\infty}e^{-s\tau}\psi_{d}(\theta,\tau)\,d\tau 
             \end{equation}
             \begin{equation}
             \hat{\rho}_{d}(\theta,s)\rule{1mm}{0mm} = \int_{0}^{\,\infty}e^{-s\tau}\rho_{d}(\theta,\tau)\,d\tau \nonumber  \,,
             \end{equation}
             and with transform variable $s$ initially required to have its real part positive.  Indeed, we see at
             once, with $\tau$ set equal to zero, theat (4) may now be simply rephrased as
             \begin{equation}          
             \psi_{d}(\theta,0)=-\frac{\omega}{2\pi\mu}\,\hat{\rho}_{d}(-1/\mu)+\frac{\omega}{2\pi}
             \frac{\mu_{0}}{\mu_{0}-\mu} \,.
             \end{equation}             
             Laplace transformation of (3) accordingly gives
             \begin{equation}
             -\mu\psi_{d}(\theta,0)+\left\{\rule{0mm}{3.5mm}1+s\mu\right\}\hat{\psi}_{d}(\theta,s)=
             \frac{\omega}{2\pi}\,\hat{\rho}_{d}(s)+\frac{\omega}{2\pi}\frac{\mu_{0}}{1+s\mu_{0}}\,,
             \end{equation}
             which can further be integrated over angle so as to yield
             \begin{equation}
             \left\{\rule{0mm}{5mm}1-\frac{\omega}{2\pi}\int_{-\pi}^{\,\pi}\frac{d\theta}{1+s\mu}\right\}\hat{\rho}_{d}(s)=
             \frac{\omega}{2\pi}\frac{\mu_{0}}{1+s\mu_{0}}\int_{-\pi}^{\,\pi}\frac{d\theta}{1+s\mu}+
             \int_{\pi/2<|\theta|<\pi}\frac{\,\mu\psi_{d}(\theta,0)\,}{\,1+s\mu\,}d\theta\,.
             \end{equation}
             The second integral on the right incorporates the null re\"{e}ntrant flux condition at the half-space
             boundary $\tau=0.$  Since it requires knowledge about the albedo, it is {\em{a priori}} unknown.
             Such information lapse notwithstanding, it is evident by inspection that that second integral is analytic
             in a left half-plane having $\Re\, s < 1.$
             
                 We embark next on a standard Wiener-Hopf journey wherein the ingredients of (9) are suitably
                 rearranged as to their left/right half-planes of analyticity, the end result being a recognition
                 of (9) as a bridge uniting ostensibly disparate, left/right analytic functions into one global
                 entire function bounded at infinity and therefore, by virtue of Liouville's theorem, a mere constant.
                 Self-evident asymptotic estimates require that constant to be null, at which point all solution
                 details fall into their rightful place.
                 
                 We begin with the definite integral
                 \begin{equation}
                 I=\int_{-\pi}^{\,\pi}\frac{d\theta}{1+s\mu}\,,
                 \end{equation}
                 which is clearly analytic when $|\Re\, s|<1,$ which we now assume, but, strictly speaking,
                 is undefined without further qualification when instead $\Im\,s=0$ whereas $|\Re\,s>1.$
                 It can be rewritten in the usual way
                 \begin{equation}
                 I=\frac{2}{\,is\,}\oint_{|\zeta|=1}\frac{d\zeta}{\zeta^{2}+2\zeta/s+1}
                 \end{equation}                                
                 as a contour integral around the unit disk $|\zeta|=|e^{i\theta}|=1$ and thus evaluated by
                 adding residues at one, two, or none of its simple poles at
                 \begin{equation}
                 \zeta_{\pm}=\left\{\rule{0mm}{3mm}-1\pm i\sqrt{s^{2}-1\,}\right\}/s\,.
                 \end{equation}
                 \newpage
                 \mbox{    }
                 \newline
                 \newline
                 \newline
                 \parindent=0 in 
                 In keeping with the precautionary remarks above, we provide $\sqrt{\zeta^{2}-1\,}=\sqrt{s-1\,}\sqrt{s+1\,}$
                 with branch cuts radiating outward from $\pm1$ to $\pm\infty.$  In particular, $i\sqrt{s^{2}-1\,}$ is real
                 and negative when $\Im\,s = 0$ whereas $-1<\Re\,s<1.$  Since also $\zeta_{+}\zeta_{-}=1,$
                 it follows that only the one pole at $\zeta_{-}$ is actually enclosed.  And so
                 \begin{equation}
                 I=\frac{2}{\,is\,}\oint_{|\zeta|=1}\frac{d\zeta}{\,(\zeta-\zeta_{+})(\zeta-\zeta_{-})\,}=
                 \frac{4\pi}{\,s(\zeta_{-}-\zeta_{+})}=\frac{2\pi i}{\,\sqrt{s^{2}-1\,}\,}\,.
                 \end{equation}
                 \parindent=0.5 in
                 
                 Returning now to the main theme, and with the information from (13) explicitly displayed, we execute our first
                 rearrangement by adding and subtracting $1/(1+s\mu_{0})$ on the right in (9) so as to arrive at
                 \begin{equation}
                 \left\{\rule{0mm}{5mm}\!\frac{\,\sqrt{s^{2}-1\,}-i\omega\,}{\,\sqrt{s^{2}-1\,}\,}\!\right\}
                 \tilde{\rho}_{d}(s)=g(s)
                 \end{equation}
                 wherein we further abbreviate by setting
                 \begin{equation}
                 \tilde{\rho}_{d}(s)=\hat{\rho}_{d}(s)+\frac{\mu_{0}}{1+s\mu_{0}}\,,
                 \end{equation}
                 \begin{equation}
                 g(s)=\frac{\mu_{0}}{1+s\mu_{0}}+
                 \int_{\pi/2<|\theta|<\pi}\frac{\,\mu\psi_{d}(\theta,0)\,}{\,1+s\mu\,}d\theta\,.	
                 \end{equation}
                 As is easily ascertained, the numerator on the left in (14) has two simple zeros at\footnote{In an
                 appendix we show that $s_{\pm}$ are the propagation constants associated with source-free
                 solutions in a medium which need not necessarily occupy all of space.  In particular, $s_{-},$
                 conveying an unbounded spatial growth $\exp(-s_{-}\tau)$ for $\tau\rightarrow +\infty,$ has
                 relevance for the kindred Milne problem, wherein it represents some unspecified, deeply buried
                 photon source of unlimited strength.}
                 \begin{equation}
                  s_{\pm}=\pm\sqrt{1-\omega^{2}\,}
                 \end{equation}
                 whose presence dampens the prospect of its imminent use as the argument of a logarithm.
                 We correct for                                            
                 this by first dividing out $(s-s_{+})(s-s_{-})=s^{2}-s_{+}^{2}$ and then compensating
                 for its adverse effect as $s\rightarrow\infty$ by multiplying with
                 $s^{2}-1.$\footnote{Beyond this restricted purpose, multiplier $s^{2}-1$ evidently serves to
                 suppress denominator infinities at branch points $s=\pm1.$}  And so we mold (14) into
                \begin{equation}
                \left(\!\frac{\,s^{2}-s_{+}^{2}\,}{\,s^{2}-1\,}\!\right)\kappa(s)\,\tilde{\rho}_{d}(s)=g(s)
                \end{equation}
                with
                \begin{equation}
                \kappa(s)=\frac{\,(s^{2}-1)\left\{\rule{0mm}{3mm}\sqrt{s^{2}-1\,}-i\omega\right\}\,}
                {\,(s^{2}-s_{+}^{2})\sqrt{s^{2}-1\,}\,}\,.
                \end{equation}
                \vspace{-3mm}
                
                One immediate, but only partial step toward rearranging (18) so as to achieve a left-hand side analytic in a
                right half-plane, and a right-hand side analytic in a left plane is to set
                \begin{equation}
                \left(\!\frac{\,\,s+s_{+}\,}{\,s+1\,\,\,}\!\right)\kappa(s)\,\tilde{\rho}_{d}(s)=
                \left(\!\frac{\,s-1\,}{\,\,\,\,s-s_{+}\,}\!\right)g(s)\,.
                \end{equation}
                The real difficulty lies in achieving a similar split for function $\kappa(s),$ a task
                accomplished via the known recipe of passing first, seemingly redundantly, to the exponential of
                its logarithm, and then representing that latter as a Cauchy contour integral with up/down
                legs of infinite extent parallel to the imaginary axis of $s.$  And so, if we set, for some
                real positive $\beta,$ $0<\beta<s_{+},$\footnote{The geometric layout of quantities $\beta$
                and $s_{+}$ is illustrated in Figure 2.}	
                \begin{equation}
                \kappa_{\pm}(s)=\exp\left(\frac{1}{\,2\pi i\,}\int_{\,\pm\beta-i\infty}^{\,\pm\beta+i\infty}
                \frac{\,\rm{Log}\,\kappa(\zeta)\,}{\,\zeta-s\,}\,d\zeta\right)\,,
                \end{equation}
                \newpage
                \mbox{    }
                \newline
                \newline
                \newline
                \parindent=0 in
                then we obtain at once the decomposition
                \vspace{-3mm}
                \begin{equation}
                \kappa(s)=\frac{\,\kappa_{+}(s)\,}{\,\kappa_{-}(s)\,}\,,
                \end{equation}                              
                valid at least in the vertical strip $-\beta<\Re\, s<\beta.$  Moreover, functions $\kappa_{\pm}(s)$
                are readily seen to be analytic respectively across the overlapping left/right half planes, and
                indeed to vanish on approach to infinity, $|s|\rightarrow\infty,$ $\Re\, s < \beta$ or else
                $\Re\, s > -\beta,$ in that order.  Equality (20), when rewritten now as
                \begin{equation}
                \left(\!\frac{\,\,s+s_{+}\,}{\,s+1\,\,\,}\!\right)\!\frac{\,\tilde{\rho}_{d}(s)\,}{\,\kappa_{-}(s)\,}=
                \left(\!\frac{\,s-1\,}{\,\,\,\,s-s_{+}\,}\!\right)\!\frac{\,g(s)\,}{\,\kappa_{+}(s)\,}\,,
                \end{equation}
                represents thus yet another step in segregating quantities that are equal within a band of overlap but,
                otherwise, assert their analyticity throughout dissimilar half planes.
                \parindent=0.5 in
                
                       Only one further adjustment remains to make that segregation complete, and that is subtraction
                from $g(s)$ of the simple pole at $s=-1/\mu_{0}.$  The subtraction must of course occur on both sides
                of (20), which thus yields
                \begin{eqnarray}
                \left(\!\frac{\,\,s+s_{+}\,}{\,s+1\,\,\,}\!\right)\!\frac{\,\tilde{\rho}_{d}(s)\,}{\,\kappa_{-}(s)\,}-
                \frac{\,\mu_{0}(1+\mu_{0})\,}{\,\kappa_{+}(-1/\mu_{0})(1+\mu_{0}s_{+})(1+\mu_{0}s)\,} &  &\nonumber\\
                     &\rule{-10cm}{0mm}=&\rule{-5cm}{0mm}\left(\!\frac{\,s-1\,}{\,\,\,\,s-s_{+}\,}\!\right)
                \!\frac{\,g(s)\,}{\,\kappa_{+}(s)\,}
                -\frac{\,\mu_{0}(1+\mu_{0})\,}{\,\kappa_{+}(-1/\mu_{0})(1+\mu_{0}s_{+})(1+\mu_{0}s)\,}
                \end{eqnarray}                               
                as a {\em{bona fide}} bounded entire function, null at infinity, and thus null everywhere on the
                strength of Liouville's theorem.  But that means that our solution is in hand, inasmuch as (15)
                gives, first                 
                \begin{equation}
                \hat{\rho}(s)=\frac{\,\mu_{0}(1+\mu_{0})\,}{\,\kappa_{+}(-1/\mu_{0})(1+\mu_{0}s_{+})\,}
                \frac{\,(s+1)\kappa_{-}(s)\,}{\,(s+s_{+})(1+\mu_{0}s)\,}-\frac{\,\mu_{0}\,}{\,(1+\mu_{0}s)\,} \,,
                \end{equation}
                and then the desired albedo follows from (7) as
                \begin{eqnarray}
                \psi_{d}(\theta,0)  &  =  & -\frac{\omega}{2\pi\mu}\,\hat{\rho}_{d}(-1/\mu)+\frac{\omega}{2\pi}
                \frac{\mu_{0}}{\mu_{0}-\mu}      \nonumber  \\
                &  =  &\frac{\omega}{2\pi}\frac{\,\mu_{0}(1+\mu_{0})(1-\mu)\kappa_{-}(-1/\mu)\,}
                {\,(\mu_{0}-\mu)(1+\mu_{0}s_{+})(1-\mu s_{+})\kappa_{+}(-1/\mu_{0})\,} \,.                 
                \end{eqnarray}
                Somewhat greater symmetry in $\mu_{0}$ and $\mu$ is attained by noting from (21) that
                \begin{equation}
                \kappa_{-}(-1/\mu)=1/\kappa_{+}(1/\mu)
                \end{equation}                              
                since, in the present instance, $\mu$ is negative and, by virtue of (19),
                $\kappa(s)$ {\em{per se}} is symmetric, $\kappa(s)=\kappa(-s),$ under argument reflection
                through the origin.  On replacing $\mu$ with $-\mu_{0}$ (27) also gives
                \begin{equation}
                \kappa_{-}(1/\mu_{0})=1/\kappa_{+}(-1/\mu_{0})\,.
                \end{equation}
                Hence
                \begin{eqnarray}
                \psi_{d}(\theta,0) & = &\frac{\omega}{2\pi}\frac{\,\mu_{0}(1+\mu_{0})(1-\mu)\,}
                {\,(\mu_{0}-\mu)(1+\mu_{0}s_{+})(1-\mu s_{+})\,}\times \kappa_{-}(1/\mu_{0})\kappa_{-}(-1/\mu) \nonumber  \\
                                   &   &                                                            \\
                                   & = &\frac{\omega}{2\pi}\frac{\,\mu_{0}(1+\mu_{0})(1-\mu)\,}
                {\,(\mu_{0}-\mu)(1+\mu_{0}s_{+})(1-\mu s_{+})\,}\times\frac{\,1\,}{\,\kappa_{+}(-1/\mu_{0})\kappa_{+}(1/\mu)\,}\,.   \nonumber
                \end{eqnarray}
                \newpage
                \mbox{   }
                \newline
                \newline
                \newline                
                \parindent=0 in
                Furthermore, if we set\footnote{In a perhaps somewhat nonstandard notation, $D$ stands for {\em{diffuse.}}}
                \begin{equation}
                \psi_{d}(\theta,0)=-\mu^{-1}D(-\mu,\mu_{0})=-\frac{\,1\,}{\,\mu\,}\left\{\rule{0mm}{6mm}\!
                -\frac{\omega}{2\pi}\frac{\,\mu_{0}\mu(1+\mu_{0})(1-\mu)\,}
                {\,(\mu_{0}-\mu)(1+\mu_{0}s_{+})(1-\mu s_{+})\kappa_{+}(-1/\mu_{0})\kappa_{+}(1/\mu)\,}\,\right\}\,,
                \end{equation}
                then it follows at once that
                \begin{equation}
                D(\mu_{0},-\mu)=D(-\mu,\mu_{0})\,,
                \end{equation}
                a symmetry under interchange of incoming/outgoing photon directions at the half-space boundary,
                traditionally associated with invariant embedding [{\textbf{7-10}}].  With a view to the first line
                in (29) it remains now only to
                deform the vertical contour whereby $\kappa_{-}(s)$ is defined in (21) into a form more
                convenient for numerical implementation.
                \parindent=0.5 in
                \subsection{Contour deformation}
                
                       Of the two options among $\kappa_{\pm}$ which Eq. (29) provides, we choose the first
                       and proceed to deform its contour around the branch cut extending from $s=-1$ to
                       $s=-\infty,$ as shown schematically on Figure 2.  Potential contributions from both
                       the semi-circle at infinity and from the full, but vanishingly small circle around the branch point at $s=-1$ are dismissed in the usual way, and what remains is the difference of integrating
                       immediately above and below the branch cut.\footnote{One must resist the natural temptation of confusing the half-plane analyticity regions for $\kappa_{\pm}(s)$ with the more restricted behavior of $\kappa(s)$ {\em{per se}}.  Indeed, save for their individual, vertical integration paths, each of $\kappa_{\pm(s)}$ defines {\em{two}} analytic functions in adjacent half planes filling up the entire two dimensional space.  But these are niceties without physical relevance to our problem, and will thus not be discussed.}  We thus find in short order that                 
                 \vspace{-0.25in}
                 \begin{center}
                 	\includegraphics[angle=-90,width=0.65\linewidth]{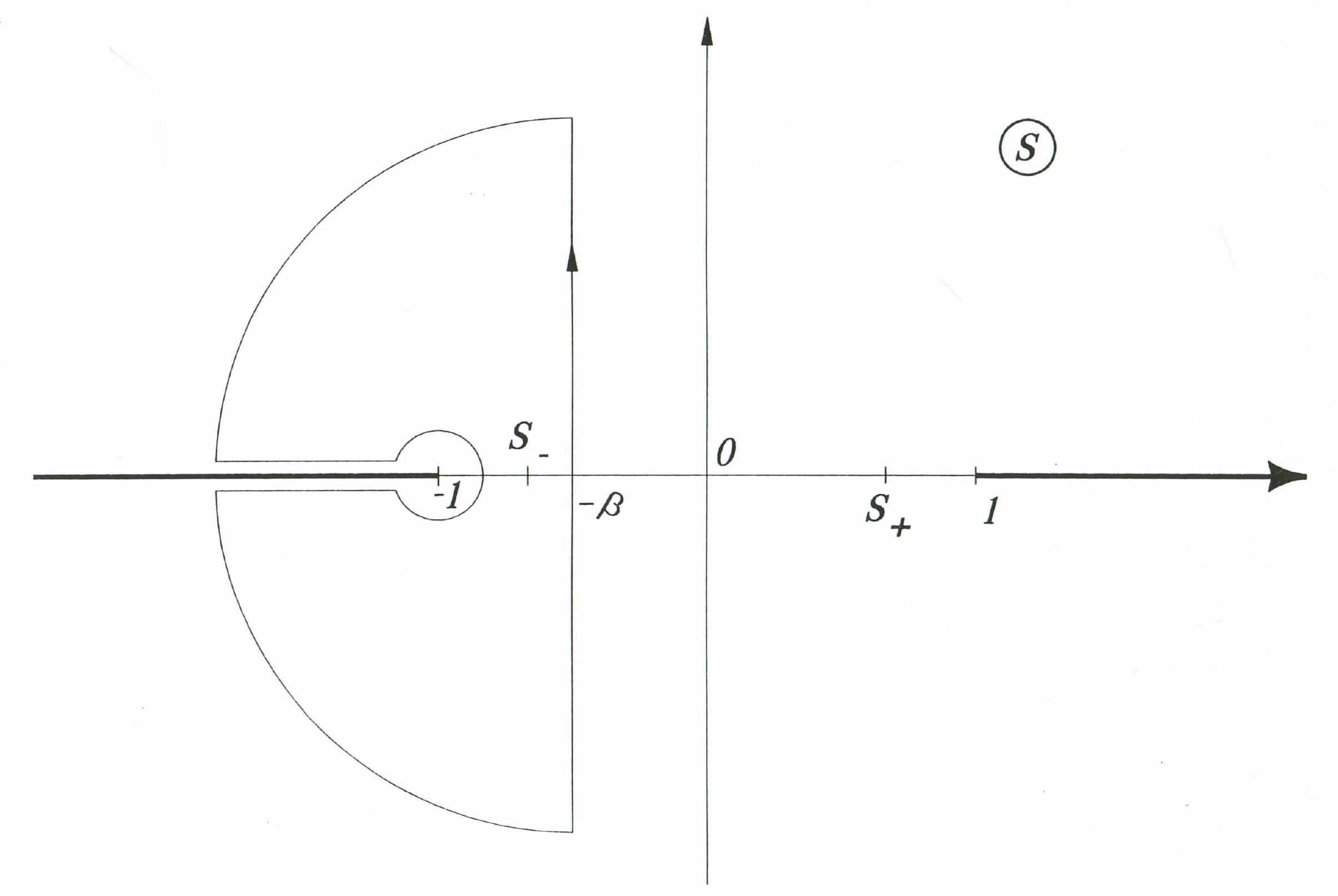}
                 \end{center}
                 \vspace{-7mm}
                 \begin{center}
                 	{\large{
                 			Figure 2.  Contour deformation for calculating $\kappa_{-}(s)$ } }
                 \end{center}
                 \vspace{-3mm}
                 \begin{eqnarray}
                 \kappa_{-}(s) & = & \exp\left[\rule{0mm}{7mm}\frac{\,1\,}{\,\pi\,}\int_{\,1}^{\,\infty}\!
                                   \tan^{-1}\left\{\!\frac{\,\omega\,}{\,\sqrt{\rule{0mm}{3.15mm}\zeta^{2}-1\,}\,}\!\right\}
                                   \!\frac{\,d\zeta\,}{\,\zeta+s\,}\right]\nonumber\\
                             & = & \exp\left[\rule{0mm}{7mm}\frac{\,\omega^{2}\,}{\,\pi\,}\int_{\,0}^{\,\pi/2}\!\!
                             \frac{\,\cot\vartheta\csc\vartheta\,}
                                  {\,\left(s+\sqrt{\rule{0mm}{3.5mm}1+\omega^{2}\cot^{2}\vartheta\,}\,\right)\!
  	                              \sqrt{\rule{0mm}{3.25mm}1+\omega^{2}\cot^{2}\vartheta\,}\,}\,       
                                  \vartheta d\vartheta\right]        \,,    
                 \end{eqnarray}
                 \newpage
                 \mbox{   }
                 \newline
                 \newline
                 \newline                
                 \parindent=0 in
                 the second line following under the natural variable substitution
                 \[\frac{\,\omega\,}{\,\sqrt{\rule{0mm}{3.15mm}\zeta^{2}-1\,}\,}=\tan\vartheta\,.\]
                 In its latter guise the computation of $\kappa_{-}(s)$ becomes amenable to numerical integration.
                 \parindent=0.5 in
                 
                 Entirely similar reasoning, but with the contour wrapped now around the complementary branch cut
                 extending from $s=1$ to $s=\infty,$ yields                       
                 \begin{eqnarray}
                 \kappa_{+}(s) & = & \exp\left[\rule{0mm}{7mm}-\frac{\,1\,}{\,\pi\,}\int_{\,1}^{\,\infty}\!
                 \tan^{-1}\left\{\!\frac{\,\omega\,}{\,\sqrt{\rule{0mm}{3.15mm}\zeta^{2}-1\,}\,}\!\right\}
                 \!\frac{\,d\zeta\,}{\,\zeta-s\,}\right]    \nonumber\\
                 & = & \exp\left[\rule{0mm}{7mm}\frac{\,\omega^{2}\,}{\,\pi\,}\int_{\,0}^{\,\pi/2}\!\!
                 \frac{\,\cot\vartheta\csc\vartheta\,}
                 {\,\left(s-\sqrt{\rule{0mm}{3.5mm}1+\omega^{2}\cot^{2}\vartheta\,}\,\right)\!
                 	\sqrt{\rule{0mm}{3.5mm}1+\omega^{2}\cot^{2}\vartheta\,}\,}\,       
                 \vartheta d\vartheta\right]        \,.    
                 \end{eqnarray}
                 And then, if we recall that (32) holds for $\Re s> -\beta$ whereas (33) requires instead that
                 $\Re s < \beta,$ it becomes evident that Eqs. (32)-(33) abide by the symmetry embodied in
                 Eqs. (27)-(28) and, indeed, can be gotten on the basis of that symmetry, one from the other,
                 without any additional calculations.\footnote{There is a considerable amount of hidden latitude
                 in the choice of vertical contours in (21).  Within obvious limits they need not be placed at
                 strictly the same distance $\beta,$ left and right, from the imaginary axis, and, again within
                 limits, they could be allowed to undulate to a certain extent while maintaining overall a steady
                 upward progress.  But these are totally dispensable frivolities, having impact neither upon
                 the transport phenomenology, nor upon the mathematics.}
                                
              \subsection{Appendix}
              \subsubsection{Free modes}
               In Eq. (17) and its accompanying Footnote 8 we had already claimed that quantities
               $s_{\pm}=\pm\sqrt{1-\omega^{2}\,}$ are associated with source-free transport modes.  Here
               we verify this assertion in two ways, the second of which will anticipate a Milne-type integral equation for
               photon density $\rho_{d}(\tau).$  Hence if we begin with
               \begin{equation}
               \mu\frac{\partial \psi(\theta,\tau)}{\partial \tau}+\psi(\theta,\tau)=\frac{\omega}{2\pi}
               \int_{-\pi}^{\,\pi}\psi(\theta',\tau)d\theta'
               \end{equation}
               and seek a solution having an exponential behavior $\psi(\theta,\tau)=e^{-\lambda\tau}\eta(\theta),$
               with $\lambda$ strictly real, we find in the usual way that (34) requires
               \begin{equation}
               1=\frac{\,\omega\,}{\,2\pi\,}\int_{-\pi}^{\,\pi}\frac{\,d\theta\,}{\,1-\lambda\mu\,}
               \end{equation}
               once
               \[\int_{-\pi}^{\,\pi}\eta(\theta)d\theta\]
               has been cancelled from both sides.  Reference to Eqs. (10) and (13) and to the intervening discussion shows
               next that, in order to produce on the right in (35) a result that is manifestly real, we must trap $\lambda$
               between $-1$ and $1,$ obtaining
               \begin{equation}
               1=\frac{\,\omega\,}{\,\sqrt{1-\lambda^{2}\,}\,}
               \end{equation}
               or else $\lambda=s_{\pm}$ as claimed.
               \newpage
               \mbox{   }
               \newline                           
                              
               The end goal of (17) can likewise be reached by converting (34) to an integral equation for $\rho(\tau),$
               at least when the scattering medium fills all of space, $-\infty<\tau<\infty.$
               Such conversion is easily attained by use of integrating factors, which give
               \begin{equation}
               \psi(\theta,\tau)=\frac{\,\omega\,}{\,2\pi\mu\,}\int_{-\infty}^{\,\tau}\rho(\tau')e^{-(\tau-\tau')/\mu}d\tau'
               \end{equation} 
               when $\mu>0,$ and
               \begin{equation}
               \psi(\theta,\tau)=\frac{\,\omega\,}{\,2\pi|\mu|\,}\int_{\,\tau}^{\,\infty}\rho(\tau')e^{-(\tau'-\tau)/|\mu|}d\tau'
               \end{equation}
               if instead $\mu<0,$ the distinction embodying a self-evident physical causality.  Integrating over all angles we thus obtain
               \begin{equation}
               \rho(\tau)=\frac{\,\omega\,}{\,2\pi\,}\int_{\,-\infty}^{\,\infty}\rho(\tau')
               \left\{\rule{0mm}{6mm}\int_{-\pi/2}^{\,\pi/2}e^{-|\tau'-\tau|/\mu}\mu^{-1}d\theta\right\}d\tau'\,.
               \end{equation}
               Setting $\rho(\tau)=\alpha\, e^{-\lambda\tau},$ with $\alpha$ some irrelevant but positive constant, and
               real $\lambda$ less than $1$ in magnitude, $-1<\lambda<1,$ so as to assure quadrature convergence, we
               find next, with interchange of integration order taken for granted,\footnote{That is to say, Fubini's theorem
               	is assumed to be in force.}
               \begin{equation}
               \int_{-\infty}^{\,\tau}e^{-\lambda\tau'-(\tau-\tau')/\mu}d\tau'=\frac{\,\mu\,}{\,1-\lambda\mu\,}\,e^{-\lambda\tau}
               \end{equation}               
               and
               \begin{equation}
               \int_{\,\tau}^{\,\infty}e^{-\lambda\tau'-(\tau'-\tau)/\mu}d\tau'=
               \frac{\,\mu\,}{\,1+\lambda\mu\,}\,e^{-\lambda\tau}\,.
               \end{equation}
               And then, with (40)-(41) brought to bear upon (39) we find               
               \begin{eqnarray}
               1 & = & \frac{\,\omega\,}{\pi\,}\int_{-\pi/2}^{\,\pi/2}\frac{\,d\theta\,}{\,1-\lambda^{2}\mu^{2}\,} \nonumber\\
                 & = & \frac{\,\omega\,}{\,\pi(2-\lambda^{2})\,}\int_{-\pi}^{\,\pi}
                 \frac{\,d\theta\,}{\,1-\lambda^{2}\mu/(\,2-\lambda^{2}\,)\,}         \\
                 & = & \frac{\,\omega\,}{\,\sqrt{\rule{0mm}{3mm}1-\lambda^{2}\,}\,}    \nonumber
               \end{eqnarray}
               which is nothing other than (36) once more.\footnote{The middle line of (42) follows from a routine appeal to
               a half-angle trigonometric formula, whereas its last utilizes the evaluation of (10) as provided by (13).}
            \subsubsection{Milne-type integral equation}
            
            Equations (37)-(39) suggest that an integral equation can likewise be contrived for the somewhat more physically relevant
            diffuse density $\rho_{d}(\tau),$ associated with Eq. (3) and having only the half-line $\tau\geq0$ for its support.
            Analogues to (37)-(38) emerge now as
            \begin{equation}
            \psi_{d}(\theta,\tau)=\frac{\,\omega\,}{\,2\pi\mu\,}\int_{\,0}^{\,\tau}\rho(\tau')e^{-(\tau-\tau')/\mu}d\tau'+
            \frac{\,\omega\,}{\,2\pi\,}\frac{\,\mu_{0}\,}{\,\mu_{0}-\mu\,}
            \left(\rule{0mm}{5mm}e^{-\tau/\mu_{0}}-e^{-\tau/\mu}\right)
            \end{equation} 
            when $\mu>0,$ and
            \begin{equation}
            \psi_{d}(\theta,\tau)=\frac{\,\omega\,}{\,2\pi|\mu|\,}
            \int_{\,\tau}^{\,\infty}\rho_{d}(\tau')e^{-(\tau'-\tau)/|\mu|}d\tau'+
            \frac{\,\omega\,}{\,2\pi\,}\frac{\,\mu_{0}\,}{\,\mu_{0}-\mu\,}e^{-\tau/\mu_{0}}
            \end{equation}
            \newpage
            \mbox{   }
            \newline
            \newline
            \newline                
            \parindent=0 in
            when $\mu<0.$  Integration over all angles then gives
            \begin{eqnarray}
            \rho_{d}(\tau)& = &\frac{\,\omega\,}{\,\pi\,}\int_{\,0}^{\,\infty}\rho_{d}(\tau') 
    \left\{\rule{0mm}{6mm}\int_{\,0}^{\,\pi/2}e^{-|\tau'-\tau|/\mu}\mu^{-1}d\theta\right\}d\tau'   \nonumber \\
     &    &    +\frac{\,\omega\mu_{0}\,}{\,\pi\,}e^{-\tau/\mu_{0}}\int_{\,0}^{\,\pi/2}
            \frac{\,d\theta\,}{\,\mu_{0}+\mu\,}+
            \frac{\,\omega\mu_{0}\,}{\,\pi\,}\int_{\,0}^{\,\pi/2}
            \left(\rule{0mm}{5mm}e^{-\tau/\mu_{0}}-e^{-\tau/\mu}\right)\frac{\,d\theta\,}{\,\mu_{0}-\mu\,}\,.
            \end{eqnarray}
            Of the three angular quadratures on the right in (45), only
            \begin{equation}
            \int_{\,0}^{\,\pi/2}\frac{\,d\theta\,}{\,\mu_{0}+\mu\,}=
            2\,\frac{\,\tanh^{-1}\!\left(\sqrt{\frac{\,1-\mu_{0}\,}{\,1+\mu_{0}\,}}\,\right)\,}{\,\sqrt{1-\mu^{2}_{0}\,}\,}=
            \csc|\theta_{0}|\ln\left\{\rule{0mm}{6mm}\frac{1+\tan(|\theta_{0}|/2)}{1-\tan(|\theta_{0}|/2)}\right\}
            \end{equation}
            submits to a closed-form evaluation,\footnote{Courtesy of the Wolfram {\em{Mathematica}} Integrator, available
            {\em{gratis}} online.} the other
            two\footnote{About the third integral we can at least confidently assert,
            if nothing else, that it is positive, so that it poses no danger whatsoever of pushing $\rho_{d}(\tau)$ into negative
            territory.  Of course, both first and second angular quadratures are positive by inspection.}
            being considerably more recondite.  In particular, the kernel
            \begin{equation}
            K(\tau',\tau)=\int_{\,0}^{\,\pi/2}e^{-|\tau'-\tau|/\mu}\mu^{-1}d\theta=\int_{\,1}^{\infty}
            \frac{\,e^{-\zeta|\tau'-\tau|}\,}{\,\sqrt{\rule{0mm}{3.10mm}\zeta^{2}-1\,}\,}\,d\zeta
            \end{equation}            
            of integral equation (45) mimics its exponential integral analogue
            \begin{equation}
            E_{1}(\tau',\tau)=\int_{\,0}^{\infty}e^{-\zeta|\tau'-\tau|}\zeta^{-1}d\zeta
            \end{equation}
            as normally encountered in connection with the Milne problem in a {\em{bona fide}} three-dimensional, isotropic
            scattering context [{\bfseries{12}}].  Neither one of the kernels (47)-(48) seems to admit evaluation in closed
            form.  This impediment notwithstanding, a Wiener-Hopf attack can be successfully mounted in the presence of
            kernel (48), and so presumably it could likewise be pursued with that of (47), essentially scrolling the
            analysis somewhat in reverse.  But it is best at this point to ignore the siren call which beckons from deep
            within this detour.
            \parindent=0.5 in
            \section{Discrete ordinates}
                   The Wiener-Hopf apparatus pivots around the global attribute $\rho_{d}(\tau),$ which no longer retains
            detailed flux memory $\psi_{d}(\theta,\tau)$ along individual photon flight directions $\theta.$  Such details
            are restored, if only approximately, by discretizing the angular quadrature and then insisting that Eq. (3)
            remain valid at the angular nodes.\footnote{A preferred, easily implemented discretization is obtained by
            partitioning the full, $-\pi<\theta<\pi$ range into any desired number of contiguous subintervals and, within
            each such, placing down a Gauss-Legendre (GL) net of some moderately high order.  If $2a$ is the common
            subinterval length, then the GL weights $w$ are correspondingly scaled down, $w\rightarrow aw.$  The global
            quadrature itself, evidently, is obtained as a simple sum over subinterval quadratures.  Such concatenated
            GL discretization is easy to implement and has the merit of avoiding end-point bunching incident to use of
            just one master GL cyle of very high order ({\em{i.e.,}} for $\theta$ close to and above/below $\mp\pi$).
            Some interior bunching is of course inevitable, but in general there emerges a much more uniform composite
            GL net.}


                    We thus imagine the discretized angular index $k$ to run from $1$ to some
            chosen $N$\footnote{$\rule{0mm}{0.75in}$It is
            of considerable bookkeeping convenience here to have $N$ divisible by four, with the second and third quartiles
            alluding to forward photon flux, while the first and fourth, when examined at the boundary $\tau=0,$
            measure the angular albedo distribution.}         	
                    and organize the
            various discretized values $\psi_{d}(\theta_{k},\tau)$ into a column vector
            \begin{equation}
            \Psi_{d}(\tau)=\left[ \begin{array}{c}
                                 \psi_{d}(\theta_{1},\tau)    \\
                                 \psi_{d}(\theta_{2},\tau)    \\
                                       \vdots        \\
                                       \vdots        \\
                                \psi_{d}(\theta_{N-1},\tau)  \\
                                 \psi_{d}(\theta_{N},\tau)
                                \end{array}   \right]\,.
             \end{equation}
             \newpage
             \mbox{   }
             \newline
             \newline
             \newline                
             \parindent=0 in
             Also required for the source term in (3) is the vector
             \begin{equation}
             M          =\left[ \begin{array}{c}
                                       1/\mu_{1}       \\
                                       1/\mu_{2}       \\
                                       \vdots          \\
                                       \vdots          \\
                                       1/\mu_{N-1}     \\
                                       1/\mu_{N}
                                \end{array}   \right]
             \end{equation}             
             and an $N\times N$ matrix
             \begin{equation}
               S_{k,l}=\frac{\,\delta_{kl}-\omega aw_{l}/2\pi\,}{\,\mu_{k}\,}  
             \end{equation}
             which accounts for the competition between absorption loss and scattering gain.\footnote{In (51), 
             $\delta_{kl}$ is the familiar Kronecker delta, equal to $1$ when $k=l$ and $0$ otherwise.}
         \parindent=0.5 in
            
                   It is a matter of empirical evidence that, in the Fortran codes which have been written
             around these ideas, all eigenvalues $\lambda_{k}$ of $S,$
             \begin{equation}
             S\Lambda_{k}=\lambda_{k}\Lambda_{k}
             \end{equation}
             with $1\leq k \leq N$ and $\Lambda_{k}$ being the corresponding eigenvector, have invariably
             turned out to be both real and distinct.  We are admittedly unable to prove now that such
             a circumstance will prevail under all admissible parameter entries, but, should its truth
             be provisionally accepted as universal, then we are guaranteed by a well known theorem that
             the eigenvector set $\{\Lambda_{k}\}_{k=1}^{N}$ provides a basis for the $N-$dimensional linear
             vector space before us.  Incidentally, the eigenvalues $\{\lambda_{k}\}_{k=1}^{N}$ themselves
             have always been found to split evenly into positive and negative values symmetrically deployed
             around the origin, $\{\lambda_{k}\}_{k=1}^{N/2}$ all negative, their complement
             $\{\lambda_{k}\}_{k=N/2+1}^{N}$ all positive.  The eigenvalue/eigenbasis framework thus
             revealed dominates all developments about to ensue, and, in order to impart a concrete symbolism
             to that framework, we set
             \begin{equation}
             \Lambda_{k}         =\left[ \begin{array}{c}
             \Lambda_{1,k}       \\
             \Lambda_{2,k}       \\
             \vdots          \\
             \vdots          \\
             \Lambda_{N-1,k} \\
             \Lambda_{N,k}
             \end{array}   \right] \,.
             \end{equation}

                   We thus duly arrive at
             \begin{equation}
             \frac{\,d\,}{\,d\tau\,}\left\{e^{\tau S}\Psi_{d}(\tau)\right\}=
             \frac{\,\omega\,}{\,2\pi\,}e^{-\tau/\mu_{0}}e^{\tau S}M=\frac{\,\omega\,}{\,2\pi\,}e^{-\tau/\mu_{0}}
             \sum_{k=1}^{N}e^{\lambda_{k}\tau}m_{k}\Lambda_{k}
             \end{equation}      
             as the discretized counterpart of (3), with real numbers $\{m_{k}\}_{k=1}^{N}$ being the
             expansion coefficients of source vector (50) in our newly discovered eigenbasis, obtained by
             standard methods of linear algebra.  Integrating
             upward from $\tau=0$ gives\footnote{As was similarly found in connection with Eqs. (43) and (45),
             neither does there intrude into (54) any need to fear encroachment by some one $\lambda_{k}>0$
             upon $\mu_{0}.$  Indeed, in both cases l'H\^{o}pital's rule saves the day.  But, should such l'H\^{o}pital
             limit extraction ever be required, we would still be assured of an ultimate asymptotic decay,
             proceeding now as $\propto \tau e^{-\tau/\mu_{0}}.$} 
             \begin{equation}
             \Psi_{d}(\tau)=e^{-\tau S}\Psi_{d}(0)+
             \frac{\,\omega\mu_{0}\,}{\,2\pi\,}
             \sum_{k=1}^{N}\frac{\,m_{k}\,}{\,\lambda_{k}\mu_{0}-1\,}
             \left(e^{-\tau/\mu_{0}}-e^{-\lambda_{k}\tau}\right)\Lambda_{k}\,.
             \end{equation}
             \newpage
             \mbox{   }
             \newline
                                       
                    If we next expand $\Psi_{d}(0)$ along the eigenbasis,
             \begin{equation}
             \Psi_{d}(0)=\sum_{k=1}^{N}b_{k}\Lambda_{k}\,,
             \end{equation}              
             wherein the various amplitudes $\{b_{k}\}_{k=1}^{N}$ allude to {\em{``boundary,"}} we find that
             \begin{equation}
             \Psi_{d}(\tau)=\sum_{k=1}^{N}e^{-\lambda_{k}\tau}b_{k}\Lambda_{k}+
             \frac{\,\omega\mu_{0}\,}{\,2\pi\,}
             \sum_{k=1}^{N}\frac{\,m_{k}\,}{\,\lambda_{k}\mu_{0}-1\,}
             \left(e^{-\tau/\mu_{0}}-e^{-\lambda_{k}\tau}\right)\Lambda_{k}\,.
             \end{equation}
             Diverging exponentials with $\lambda_{k}<0$ are banished at once by setting, with
             complete impunity as to any menace of having to cope with vanishing denominators,
             \begin{equation}
             b_{k}=\frac{\,\omega\mu_{0}m_{k}\,}{\,2\pi(\lambda_{k}\mu_{0}-1)\,}
             \end{equation}
             for $1\leq k \leq N/2,$ while the remaining amplitudes $b_{k}$ at indices $N/2+1 \leq k \leq N$ follow by
             suppressing any re{\"{e}}ntrant flux,
             \begin{equation}
             \sum_{k=N/2+1}^{N}b_{k}\Lambda_{l,k}=-\sum_{k=1}^{N/2}b_{k}\Lambda_{l,k}=
             \frac{\,\omega\mu_{0}\,}{\,2\pi\,}\sum_{k=1}^{N/2}\frac{\,m_{k}\Lambda_{l,k}\,}{\,1-\lambda_{k}\mu_{0}\,}
             \end{equation}
             required now $\forall$ $N/4+1 \leq l \leq 3N/4.$  Admittedly the $N/2\times N/2$ linear system (59) is still
             in need of its own solution, but that is a relatively standard matter {\em{vis-\`{a}-vis}} the wide
             availability of matrix inversion routines.\footnote{This fresh inversion prospect should not provoke
             consternation.  Indeed, it is understood that by this point we had already confronted
             an even more daunting, $N\times N$ matrix inversion on behalf of the $m_{k}$ amplitudes first declared in
             Eq. (54).}
                        	
             \section{Numerical examples}
                         
                   Since we have no wish to inundate the prospective reader with an avalanche of numerical data,
             we provide just two albedo examples based on the developments that have been sketched, showing a
             most welcome agreement between Wiener-Hopf {\em{versus}} discrete ordinates outcomes, an agreement
             which surely bolsters confidence in both.
             
                   The above analysis has been implemented in Fortran code with both Wiener-Hopf and discrete
             ordinates viewpoints running in parallel.\footnote{In point of fact our code, at least insofar as its
            discrete ordinates part is concerned, aims at somewhat greater generality wherein the fibered, scattering
            medium has finite thickness and is sandwiched, fore and aft, by uniform, nonscattering dielectric blankets.
            All such potential generalizations, however, are simply bypassed, through suitable parameter choice,
            when undertaking the calculations now reported in Figures 3 and 4.  One may note in passing that computer
            run times underlying such calculations are, to all intents and purposes, entirely inconsequential, even
            on a personal laptop device of modest power. }  Utility subroutines for
            eigenvalue/eigenvector computation, linear algebraic system solution, and numerical integration,
            the latter on behalf of the Wiener-Hopf albedo result (26) {\em{et seq.,}} have all been drawn from the IMSLIB libraries.
          
                   The sample computations of Figures 3 and 4 have been performed on an angular grid of $40$ concatenated
                   Gaussian level-$10$ quadratures, yielding $100$ points per quadrant.  As a way to ease visualizing a
                   performance comparison, the discrete ordinates outcomes were then simply winnowed down by a factor
                   of $5$ so as to yield just $20$ point per quadrant.  It should come as no surprise that substantial
                   absorption, as conveyed by a single scattering albedo $\omega=0.60$ in Figure 3, depresses the
                   the emergent flux across the board as compared to that from an essentially dissipation-free medium
                   ($\omega=0.99$) in Figure 4.\footnote{In both figures the incidence angle $\theta_{0}$ has been
                  	set at $45^{o}.$  But, even though we do not dwell on this aspect here, the albedo patterns are only
                  	weakly dependent upon $\mu_{0}=\cos\theta_{0}.$}  The cu-
                  \newpage
                  \mbox{  }
                  \newline
                  \newline
                  \newline                      
                   mulative albedos, perpendicular to the boundary and normalized by the similarly perpendicular influx,
                   when integrated across the indicated, $200-$point angular grid, ring in at
                   $0.29253\ldots$ in the strong absorption case (Figure 3),                                               
                   and $1.14332\ldots\ldots$ when instead a pure scattering regime predominates.  We are at a loss to provide any
                   physically credible motivation for the horizon bulge in Figure 4, but of course would welcome all
                   plausible suggestions.

             \vspace{-1.0in}
             \begin{center}
             	\includegraphics[width=0.87\linewidth]{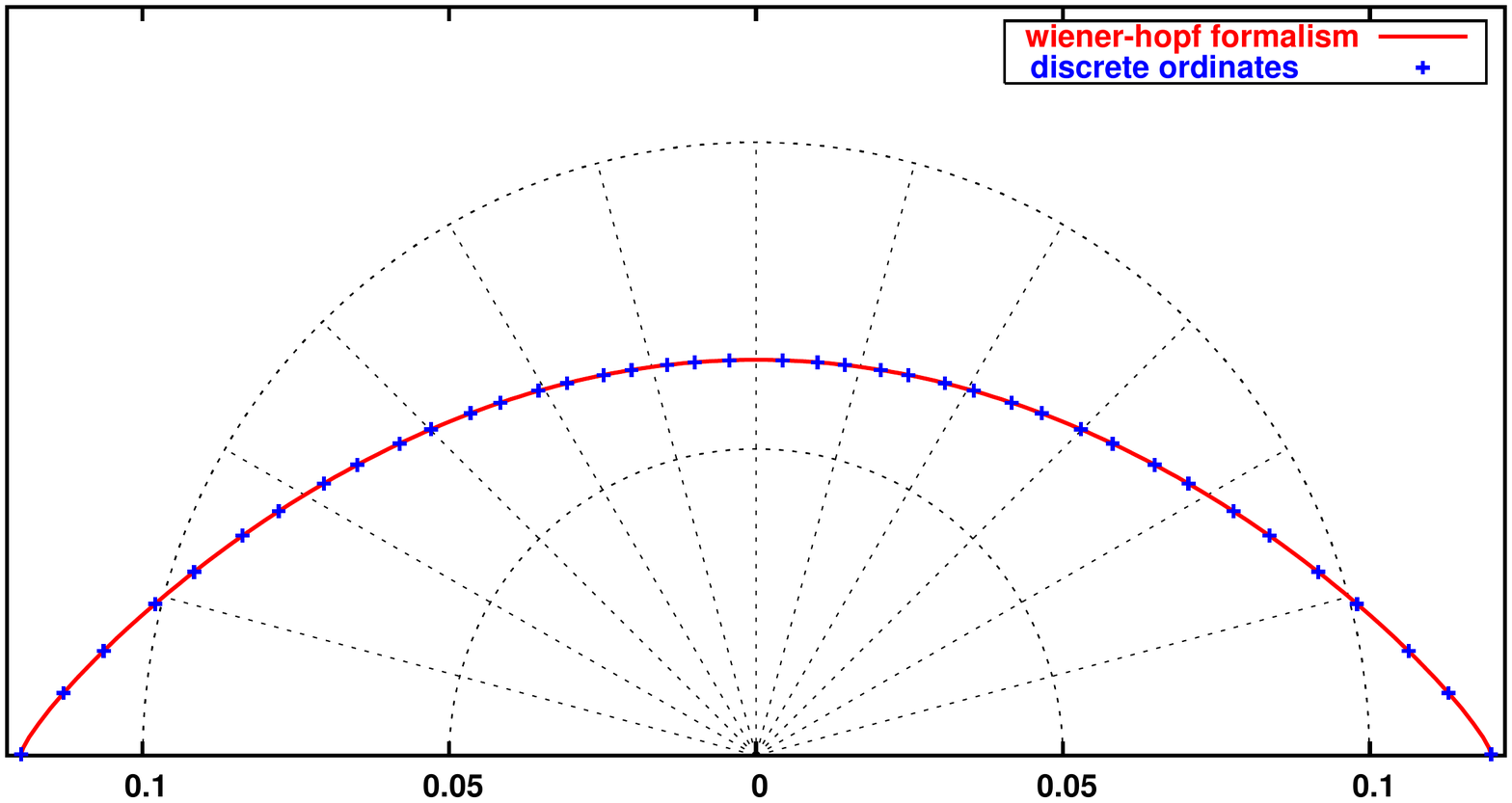}
             \end{center}
             \vspace{-35mm}
             \begin{center}
             	{\large{
             			Figure 3.  Angular albedo pattern from strongly absorbing medium \\
             			     Single scattering albedo $\omega=0.60$, incident $\theta_{0}=45$ deg \\
             			     Pattern angle $\vartheta$ measured in degrees on boundary exterior in y-x plane } }
             \end{center}
             \vspace{-1.2in}
             \begin{center}
             	\includegraphics[width=0.87\linewidth]{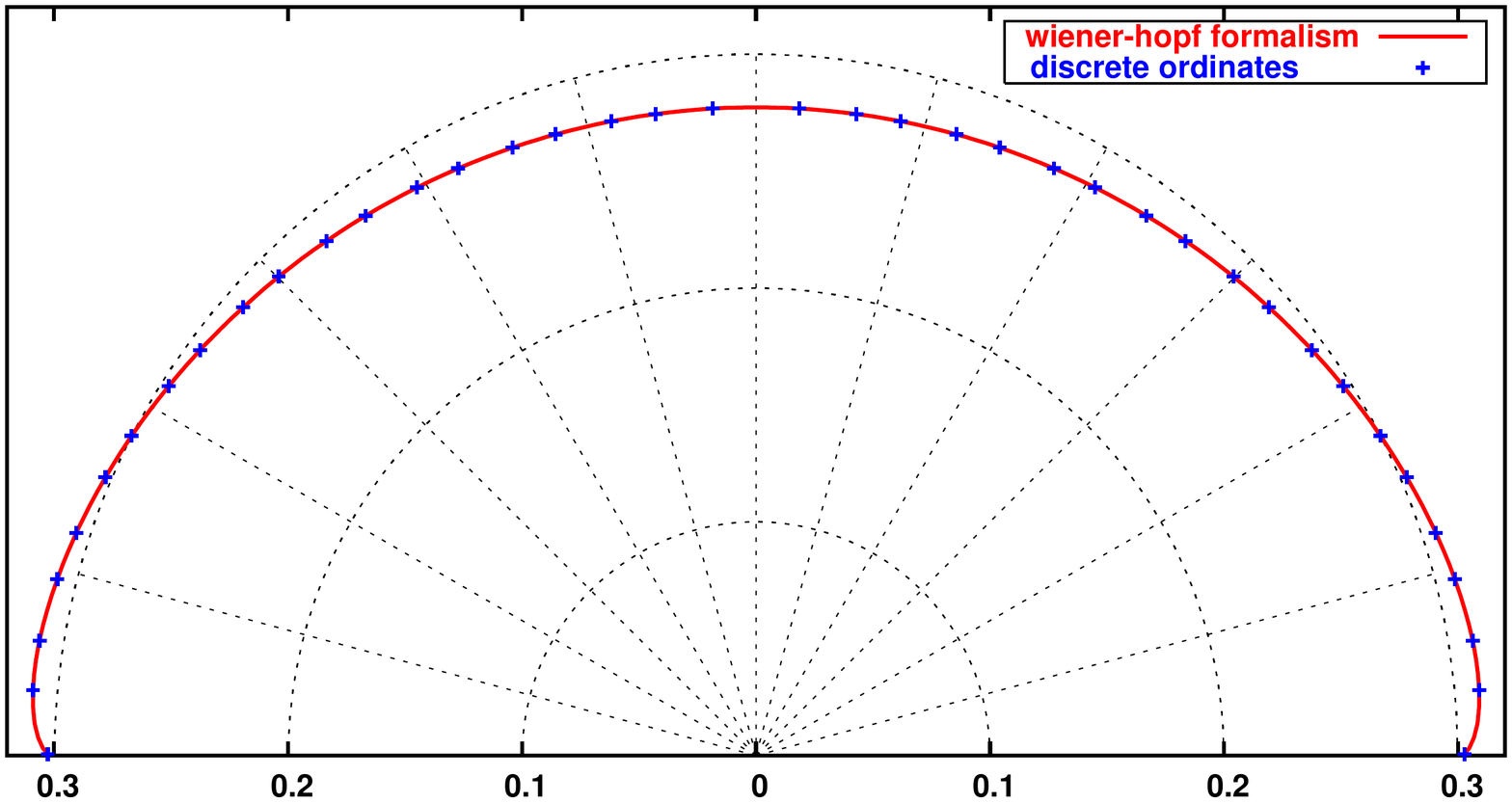}
             \end{center}
             \vspace{-35mm}
             \begin{center}
             	{\large{
             			Figure 4.  Angular albedo pattern from weakly absorbing medium \\
             			Single scattering albedo $\omega=0.99$, incident $\theta_{0}=45$ deg \\
             			Pattern angle $\vartheta$ measured in degrees on boundary exterior in y-x plane } }
             \end{center}
             \newpage
             \mbox{    }
             \section{Comments}
             \vspace{-2mm}
             
                   It is clear that the work now on view builds directly upon the standard formulation,
                   as found in, say, [{\bf{10}}] and [{\bf{12}}], wherein the scattering phenomenon is fully
                   three-dimensional, and the angular analysis is equipped with a weight appropriate to spherical
                   quadrature.  But, of course, as we all know, the devil is ever in the details, so that perhaps these
                   efforts may yet justify their existence as a model for further work by others.\vspace{-1mm}
                       
                       The next logical step in the adaptation would be to emulate a flux development along the lines
                       of Case-type singular eigenfunctions.  In fact, a first step in this direction had already been
                       taken quite a while ago, and was hastily assembled for presentation at the 1997 IEEE AP-S/URSI
                       Symposium in Montr\'{e}al, Canada [{\bf{13}}].\footnote{At first blush an IEEE/URSI conference
                       	may seem to be an unlikely, even incongruous venue for transport work of this sort.  The level
                       	of incompatibility was lessened, however, by slightly stronger allusions to the obligatory
                       	electromagnetic background, which is barely mentioned in these pages.  The present Figure 2
                       	has, incidentally, been drawn from this presentation.}  But this foray was
                       only a hint of what can and should be done, and doubtless remains to this day in need of a
                       patch up.  Perhaps such a patch up may yet be undertaken by the undersigned, and/or at the
                       hands of some other transport enthusiast(s).

                 And lastly, even though we seem to have shirked much contact with the existing transport literature,
                 we most certainly have no wish to impart the patently erroneous impression that our pages furnish
                 the one and only, the solitary foray into the realm of radiative transport across granular, inhomogeneous media.
                 On the contrary!  Our only claim is to attempt a contribution toward the understanding of radiant
                 transport wherein a prevailing background symmetry, presently of the strictly linear sort, is sufficiently
                 compelling to impress its stamp upon the propagating field structure.  Antecedents in this latter direction,
                 even if their overt concern may appear to be strictly electromagnetic, can be traced from [{\bf{4}}] and
                 [{\bf{11}}], and the present work is a direct echo thereof, transplanted into its native, transport setting.

                 By contrast, most available treatments of radiant transport across granular media confine themselves,
                 by virtue of both available analytic tools and the considerable domains of practical applicability,
                 exemplified by thermal insulation blankets and polluted atmospheres, to particle and tendril swarms,
                 the latter naturally randomized as to individual orientation.  One such study, authored by C. L. Tien,
                 a master of radiant heat transport, and T. W. Tong, one of his many students, can be found in [{\bf{14}}],
                 while a later complement, due to S. C. Lee, yet another former student, is reported as [{\bf{15}}].

                 An early treatment of radiant transport in the presence of particulate dispersal is exemplified by the 
                 Purdue University dissertation of Tom J. Love, Jr. [{\bf{16}}], wherein one finds a vast 
                 amount of numerical work implementing a variety of scattering series solutions.  The scope of this
                 numerical effort surely deserves to be judged as monumental when viewed against the relatively
                 impoverished computing resources available at the time of its creation.    By contrast, in the purely
                 theoretical work reported in[{\bf{17}}], Stephen Whitaker provides for the temperature a diffusion
                 equation whose highly modified heat flux seeks to describe porosity effects.\vspace{-1mm}

                 And, while radiant transport theory, seen from an engineering perspective, is well conveyed in
                 the standard pillars [{\bf{18}}]-[{\bf{19}}], it is only the latter which devotes much attention
                 to medium inhomogeneities.

\newpage\
                                                                                                                         
\section{References}

\parindent=0in

1.   M. Kerker, {\bf{The Scattering of Light and Other Electromagnetic Radiation}},
Academic Press, New York, 1969, pp. 255-265.

2.   M. Barabas, {\bf{Scattering of a plane wave by a radially stratified tilted cylinder}},
Journal of the Optical Society of America, Vol. A4, pp. 2240-2248 (1987).

3.   S. C. Lee, {\bf{Dependent scattering of an obliquely incident plane wave by a collection
of parallel cylinders}}, Journal of Applied Physics, Vol. 68, pp. 4952-4957 (1990).

4.   S. C. Lee, {\bf{Scattering of polarized radiation by an arbitrary collection of
closely-spaced parallel nonhomogeneous tilted cylinders}}, Journal of the Optical Society of America,
Vol. A13, pp. 2256-2265 (1996).

5.   J. A. Grzesik, {\bf{A Note on the Backward Scattering Theorem}}, Progress in
Electromagnetic Research, PIER 40, 2003, pp. 255-269.

6.   M. Kerker, D. D. Cooke, and J. M. Carlin, {\bf{Light Scattering from Infinite Cylinders.
The Dielectric-Needle Limit}}, Journal of the Optical Society of America, Vol. 60, Issue 9,
pp. 1236-1239 (1970), doi: 10.1364/JOSA.60.001236.

7. $\,$V. V. Sobolev, {\bf{A Treatise on Radiative Transfer}}, (translated by S. I. Gaposchkin),
D. Van Nostrand Company, Inc., Princeton, New Jersey, 1963.

8.  $\,\,$S. Chandrasekhar, {\bf{Radiative Transfer}}, Oxford University Press, New York, 1950.

9.  $\,\,\,$K. M. Case, {\bf{Transfer Problems and the Reciprocity Principle}}, Reviews of Modern
Physics, Vol. 29, No. 4, October, 1957, pp. 651-663.

10. $\,\,$Kenneth M. Case and Paul Frederick Zweifel, {\bf{Linear Transport Theory}}, Addison-Wesley
Publishing Company, Reading, Massachusetts, 1967.

11. $\,\,$Siu-Chun Lee and Jan A. Grzesik, {\bf{Light scattering by closely spaced parallel cylinders embedded in
a semi-infinite dielectric medium}}, Journal of the Optical Society of America, Vol. A15, pp. 163-173 (1998).

12. $\,\,$B. Davison, {\bf{Neutron Transport Theory}}, with the collaboration of J. B. Sykes, Oxford
University Press, 1958, p. 67ff.

13. $\!$J. A. Grzesik, {\bf{Two-Dimensional Radiation Transport across Fibrous Media}}, 1997 IEEE Antennas
and Propagation Society Symposium and URSI Radio Science Meeting, Montr\'{e}al, Canada, July 13-18, 1997.
\newpage
\mbox{  }
\mbox{  }
\newline
\newline
\newline
14. $\,\,$T. W. Tong and C. L. Tien, {\bf{Radiative heat transfer in fibrous insulations--Part I:  Analytical Study}},
ASME Journal of Heat Transfer, Vol. 105, No. 1, February 01, 1983, pp. 70-75.

15. $\,\,$S. C. Lee, {\bf{Radiative transfer through a fibrous medium:  allowance for fiber orientation}}, Journal of
Quantitative Spectroscopy and Radiative Transfer, Vol. 36, Issue 3, September 1986, pp. 253-263.

16. $\,\,$Tom J. Love, Jr., {\bf{An Investigation of Radiant Heat Transfer in Absorbing, Emitting, and Scattering Media}},
University of Oklahoma Research Institute, Norman, Oklahoma, January 1963.

17. $\,\,$Stephen Whitaker, {\bf{Radiant Energy Transport in Porous Media}}, Industrial and Engineering Chemical
Fundamentals, Vol. 19, May 1980, pp. 210-219.

18. $\,\,$Robert Siegel and John R. Howell, {\bf{Thermal Radiation Heat Transfer}}, {\em{Second Edition}},
Hemisphere Publishing Corporation, New York, 1981. 

19. $\,\,$Michael F. Modest, {\bf{Radiative Heat Transfer}}, {\em{Third Edition}}, Academic Press/Elsevier, New York, 2013,
{\em{cf.}} pp. 387-439.

\begingroup
\parindent 0pt
\parskip 2ex
\def\enotesize{\normalsize}                 
\theendnotes 
\endgroup

\end{document}